\documentclass{article}
\usepackage{amsfonts,amssymb,curves}

\usepackage{makeidx}     
\usepackage{graphicx}    

\newcommand{\R}{{\mathbb R}}  

\newtheorem{theorem}{Theorem}
\newtheorem{itlemma}{Lemma}[section] 
\newtheorem{itproposition}[itlemma]{Proposition}
\newenvironment{proposition}{\begin{itproposition}\rm}{\end{itproposition}}
\newenvironment{proof}{\noindent {\em Proof}.\
}{\hspace*{\fill}$\halmos$\medskip}

\newcommand{\text}[1]{\hbox{\rm \ #1\ \/}}
\newcommand{\be}[1]{\begin{equation}\label{#1}}
\newcommand{\ee}{\end{equation}}
\newcommand{\beqn}{\begin{eqnarray*}}
\newcommand{\eeqn}{\end{eqnarray*}}
\newcommand{\beq}{\begin{eqnarray}}
\newcommand{\eeq}{\end{eqnarray}}
\newcommand{\bl}[1]{\begin{lemma}\label{#1}}
\newcommand{\ble}[1]{\begin{lemmaex}\label{#1}}
\newcommand{\br}[1]{\begin{remark}\label{#1}}
\newcommand{\bt}[1]{\begin{theorem}\label{#1}}
\newcommand{\bd}[1]{\begin{definition}\label{#1}}
\newcommand{\bp}[1]{\begin{proposition}\label{#1}}
\newcommand{\bc}[1]{\begin{corollary}\label{#1}}
\newcommand{\bfact}[1]{\begin{fact}\label{#1}}
\newcommand{\ber}[1]{\begin{exercise}\label{#1}}
\newcommand{\bex}[1]{\begin{example}\label{#1}}
\newcommand{\bem}[1]{\begin{example}\label{#1}}  
\newcommand{\ec}{\mybox\end{corollary}}
\newcommand{\efact}{\mybox\end{fact}}
\newcommand{\eer}{\mybox\end{exercise}}
\newcommand{\eex}{\mybox\end{example}}
\newcommand{\eem}{\mybox\end{example}}
\newcommand{\el}{\mybox\end{lemma}}
\newcommand{\ele}{\mybox\end{lemmaex}}
\newcommand{\er}{\mybox\end{remark}}
\newcommand{\et}{\qed\end{theorem}}
\newcommand{\ed}{\mybox\end{definition}}
\newcommand{\ep}{\mybox\end{proposition}}
\newcommand{\epr}{\end{proof}}
\newcommand{\bpr}{\begin{proof}}

\newcommand{\ecs}{\end{corollary}}
\newcommand{\eers}{\end{exercise}}
\newcommand{\eexs}{\end{example}}
\newcommand{\eems}{\end{example}}
\newcommand{\els}{\end{lemma}}
\newcommand{\eles}{\end{lemmaex}}
\newcommand{\ers}{\end{remark}}
\newcommand{\ets}{\end{theorem}}
\newcommand{\eds}{\end{definition}}
\newcommand{\eps}{\end{proposition}}
\newcommand{\halmos}{\rule{1ex}{1.4ex}}
\newcommand{\mybox}{\hfill $\Box$} 

\newcommand{\ben}{\begin{enumerate}}
\newcommand{\een}{\end{enumerate}}

\newcommand{\bi}{\begin{itemize}}
\newcommand{\ei}{\end{itemize}}

\newcommand{\barM}{{\bar M}}

\begin{document}

\title{An Analysis of a Circadian Model Using The Small-Gain Approach
to Monotone Systems}
\author{David Angeli\\
Dip. Sistemi e Informatica \\
University of Florence,  50139 Firenze, Italy \\
{\tt angeli@dsi.unifi.it} \and
Eduardo D. Sontag%
\thanks{Supported in part by AFOSR Grant F49620-01-1-0063,
NIH Grants R01 GM46383 and P20 GM64375, and Aventis}
\\
Dept. of Mathematics\\
Rutgers University, NJ, USA\\
{\tt sontag@hilbert.rutgers.edu}}

\date{}
\maketitle 

\begin{abstract}

In this note, we show how certain properties of Goldbeter's 1995 model for
circadian oscillations can be proved mathematically, using techniques from the
recently developed theory of monotone systems with inputs and outputs.  
The theory establishes global asymptotic stability, and in particular
no oscillations, if the rate of transcription is somewhat
smaller than that assumed by Goldbeter.  This stability persists even under
arbitrary delays in the feedback loop.

\end{abstract}

\section{Introduction}

The molecular biology underlying the circadian rhythm in {\em Drosophila\/}
is the focus of a large amount of both experimental and theoretical work.
Goldbeter proposed a simple model for circadian oscillations
in~\cite{goldbeter95} (see also his book~\cite{Goldbeter}).
Although by now several more realistic models are available, in particular
incorporating other genes, this simpler model exhbits many realistic
features, such as a 24-hour period.
The key to the model is the inhibition of {\em per\/} gene transcription by
its protein product PER, forming an autoregulatory negative feedback loop.

In this note, we show how certain properties of the model can be proved
mathematically, using techniques from the recently developed theory of
monotone systems with inputs and outputs.  The theory establishes
global asymptotic stability, and in particular no oscillations, if the
rate of transcription is somewhat smaller than that assumed by
Goldbeter.  This stability persists even under arbitrary delays in the
negative feedback loop.
On the other hand, a larger --but still smaller than
Goldbeter's-- strength, in the presence of delays, results in oscillations.

The terminology and notations are as given in \cite{monotone1,monotone2}, and
are not repeated here.

\section{The Model}

\begin{figure}[ht]
\begin{center}
\setlength{\unitlength}{2500sp}%
\begin{picture}(8424,1647)(2389,-5848)
\thicklines
\put(3601,-4561){\vector(-1, 0){1200}}
\put(4201,-4561){\vector( 1, 0){1200}}
\put(9601,-4561){\vector( 1, 0){1200}}
\put(6001,-4486){\vector( 1, 0){1200}}
\put(7801,-4486){\vector( 1, 0){1200}}
\put(7201,-4636){\vector(-1, 0){1200}}
\put(9001,-4636){\vector(-1, 0){1200}}
\put(6001,-5761){\line(-1, 0){2100}}
\put(3901,-5761){\vector( 0, 1){825}}
\put(6601,-5836){\line( 1, 0){2775}}
\put(9376,-5836){\vector( 0, 1){900}}
\put(9226,-5011){\line( 0,-1){675}}
\put(9226,-5686){\vector(-1, 0){2625}}
\put(2851,-4336){$v_m$}
\put(3826,-4636){$M$}
\put(7351,-4636){$P_1$}
\put(5551,-4636){$P_0$}
\put(9226,-4636){$P_2$}
\put(4726,-4336){$k_s$}
\put(8326,-4336){$V_3$}
\put(8326,-4936){$V_4$}
\put(6526,-4936){$V_2$}
\put(6526,-4336){$V_1$}
\put(10126,-4336){$v_d$}
\put(3526,-5461){$v_s$}
\put(6150,-5836){$P_N$}
\put(8851,-5461){$k_1$}
\put(9526,-5461){$k_2$}
\end{picture}
\caption{Goldbeter's Model}
\label{reactions}
\end{center}
\end{figure}
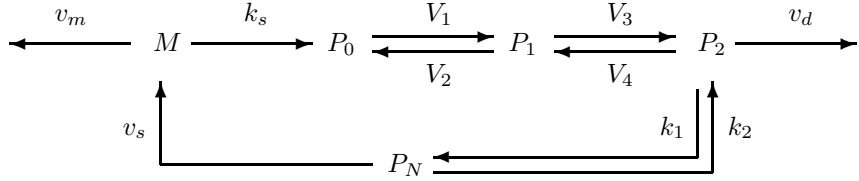

The model is as shown in Figure~\ref{reactions}.
PER protein is synthesized at a rate proportional to its mRNA concentration.
Two phosphorylation sites are available, and constitutive
phosphorylation and dephosphorylation occur with saturation dynamics,
at maximum rate $v_i$'s and with Michaelis constants $K_i$.
Doubly phosphorylated PER is degraded, also satisfying
saturation dynamics (with parameters $v_d, k_d$), and it is translocated 
to the nucleus with rate constant $k_1$.
Nuclear PER inhibits transcription of the {\em per\/} gene,
with a Hill-type reaction of cooperativity degree $n$ and threshold constant
$K_I$, and mRNA is produced. and translocated to the cytoplasm, at a rate
determined by a constant $v_s$.  Additionally, there is saturated degradation
of mRNA (constants $v_m$ and $k_m$).

The equations for concentrations are as follows:
\beqn
\dot M &=& v_sK_I^n/(K_I^n\! +\! P_N^n)-v_mM/(k_m\! +\! M)\\
\dot P_0 &=& k_sM-V_1P_0/(K_1\! +\! P_0)+V_2P_1/(K_2\! +\! P_1)\\
\dot P_1 &=& V_1P_0/(K_1\! +\! P_0)-V_2P_1/(K_2\! +\! P_1)-V_3P_1/(K_3\! +\! P_1)+V_4P_2/(K_4\! +\! P_2)\\
\dot P_2&=& V_3P_1/(K_3\! +\! P_1)-V_4P_2/(K_4\! +\! P_2)-k_1P_2+k_2P_N-v_dP_2/(k_d\! +\! P_2)\\
\dot P_N&=&k_1P_2-k_2P_N
\eeqn
where the subscript $i=0,1,2$ in the concentration $P_i$ indicates the 
degree of phosphorylation of PER protein, $P_N$ is used to indicate the
concentration of PER in the nucleus, and $M$ indicates the concentration of
{\em per\/} mRNA.
The parameters (in suitable units $\mu M$ or $h^{-1}$)
are as in the following table:
\begin{center}
\begin{tabular}{||c|c||c|c||}
\hline 
Parameter & Value & Parameter & Value \\
\hline
$k_2$ & 1.3 &
$k_1$ & 1.9 \\
$V_1$ & 3.2 &
$V_2$ & 1.58 \\ 
$V_3$ & 5 &
$V_4$ & 2.5 \\
$v_s$ & 0.76 &
$k_m$ & 0.5 \\
$k_s$ & 0.38 &
$v_d$ & 0.95 \\
$k_d$ & 0.2 &
$n$ & 4 \\
$K_1$ & 2 &
$K_2$ & 2 \\
$K_3$ & 2 &
$K_4$ & 2 \\
$K_I$ & 1 
& $v_m$ & 0.65\\
\hline
\end{tabular} 
\end{center}
With these parameters, there are limit cycle oscillations.
We leave all fixed except $v_s$, and show that there are no oscillations if
$v_s=0.4$, but oscillations exist if $v_s=0.5$ and there are delays in the
negative regulatory loop, either in transcription or in translation (or in
both). 

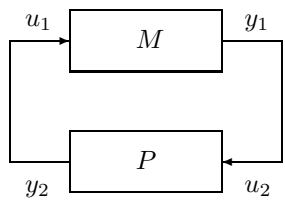
\begin{figure}[ht]
\begin{center}
\setlength{\unitlength}{2500sp}%
\begin{picture}(2724,1866)(3889,-3715)
\put(4051,-2011){$u_1$}
\put(4051,-3661){$y_2$}
\put(6226,-2011){$y_1$}
\put(6226,-3661){$u_2$}
\thinlines
\put(4501,-2461){\framebox(1500,600){}}
\put(6001,-2161){\line( 1, 0){600}}
\put(6601,-2161){\line( 0,-1){1200}}
\put(6601,-3361){\vector(-1, 0){600}}
\put(4501,-3361){\line(-1, 0){600}}
\put(3901,-3361){\line( 0, 1){1200}}
\put(3901,-2161){\vector( 1, 0){600}}
\put(4501,-3661){\framebox(1500,600){}}
\put(5150,-2236){$M$}
\put(5150,-3436){$P$}
\end{picture}
\caption{Systems in feedback}
\label{interconnection}
\end{center}
\end{figure}

We choose to view the system as the feedback interconnection
of two subsystems, see Figure~\ref{interconnection}.

\subsubsection*{mRNA System}

The first ($M$) subsystem is described by the scalar differential equation
\[
\dot M \;= \; v_sK_I^n/(K_I^n\! +\! u_1^n)-v_mM/(k_m\! +\! M)
\]
with input $u_1$ and output $y_1=k_sM$.

\subsubsection*{PER System}

The second ($P$) subsystem is four-dimensional:
\beqn
\dot P_0 &=& u_2-V_1P_0/(K_1\! +\! P_0)+V_2P_1/(K_2\! +\! P_1)\\
\dot P_1 &=& V_1P_0/(K_1\! +\! P_0)-V_2P_1/(K_2\! +\! P_1)-V_3P_1/(K_3\! +\! P_1)+V_4P_2/(K_4\! +\! P_2)\\
\dot P_2&=& V_3P_1/(K_3\! +\! P_1)-V_4P_2/(K_4\! +\! P_2)-k_1P_2+k_2P_N-v_dP_2/(k_d\! +\! P_2)\\
\dot P_N&=&k_1P_2-k_2P_N
\eeqn
with input $u_2$ and output $y_2=P_N$.

Assume from now on that:
\be{vs}
v_s \leq  0.54
\ee
(the remaining parameters will be constrained below, in such a manner that
those in the previously given table will satisfy all the constraints).

As state-space for the first system, we will pick a compact interval
$X_1=[0,\bar M]$, where
\be{ineqM}
\frac{v_sk_m}{v_m-v_s} \leq  \bar M < \frac{v_d}{k_s}
\ee
and we assume that $v_s<v_m$.
Note that the first inequality implies that
\be{eqn1}
v_s < \frac{v_m \barM}{k_m+M}
\ee
and therefore
\[
v_sK_I^n/(K_I^n\! +\! u_1^n)-v_m\barM/(k_m\! +\! \barM) < 0
\]
for all $u_1\geq 0$,
so that indeed $X_1$ is forward-invariant for the dynamics.
With the parameters shown in the table given earlier  (except for $v_s$, which
is picked as in~(\ref{vs})),
\[
\barM = 2.45
\]
satisfies all the constraints.
As input space for the mRNA system, we pick $U_1=\R_{\geq 0}$, and as output space
$Y_1=[0,v_s)$.
Note that $y_1 = k_s M \leq  k_s \barM < v_s$, by~(\ref{ineqM}), so the output
belongs to $Y_1$.

For the second system, the state space is $\R_{\geq 0}^4$, the input space
is $U_2=Y_1$, and the output space is $Y_2=U_1$.

When looking at the first system, we view $U_1$ as ordered by the cone
$\R_{\leq 0}$, but $U_2,Y_1,Y_2$ are all ordered in the usual manner
(cone $\R_{\geq 0}$).

\section{Monotonicity and Characteristics}

The first system is monotone, and has a well-defined characteristic, in the
sense of~\cite{monotone1}.  Monotonicity is clear (one-dimensional system),
and the existence of characteristics is immediate from the fact that
$\dot M>0$ for $M<k_1(u_1)$ and $\dot M<0$ for $M>k_1(u_1)$, where, for each
constant input $u_1$, 
\[
k_1(u_1)=
\frac {
v_s\,K_I^n\,k_m
}{
v_m\,K_I^n+v_m\,{u_1}^{n}-v_s\,K_I^n
}
\]
(which is an element of $X_1$).

Note that all solutions of the differential equations which describe the
$M$-system, even those that do not start in $X_1$, enter $X_1$ in finite
time (because $\dot M(t)<0$ whenever $M(t)\geq \barM$, for any input $u_1(\cdot )$).
The restriction to the state space $X_1$ (instead of using all of $\R_{\geq 0}$)
is done for convenience, so that one can view the output of the $M$ system as
in input to the $P$-subsystem.  (Desirable properties of the $P$-subsystem
depend on the restriction imposed on $U_2$.)  Given any trajectory, its
asymptotic behavior is independent on the behavior in an initial finite time
interval, so this does not change the conclusions to be drawn.  (Note that
solutions are defined for all times --no finite explosion times-- because the
right-hand sides of the equations have linear growth.)

Monotonicity of the second system is also clear, from the fact that
$\frac{\partial \dot P_i}{\partial P_j}>0$ for all $i\not= j$; in fact, this is a
{\em strongly monotone tridiagonal system\/}
(\cite{Sm, Smith}).
We show that (for the parameters in the table, as well as for a
larger set of parameters) the system has, for each constant input $u$,
a unique equilibrium, and trajectories are all bounded; it follows then
from~\cite{Sm, Smith} that the unique equilibrium is globally asymptotically
stable, which means that characteristics are well-defined.

\bp{main-P-prop}
Suppose that the following conditions hold:
\bi
\item
$v_d + V_2 < V_1$
\item
$V_1 + V_4 < V_2 + V_3$
\item
$0\leq c < v_d$
\item
$V_4+v_d<V3$
\ei
and that all constants are positive and the input $u_2(t)\equiv c$.
Then the $P$-system has a unique globally asymptotically stable equilibrium.
\eps

This will be a corollary of the following more general result.

\bt{main-P-theo}
Consider a system of the following form:
\beqn
\dot x_0&=&c-\alpha _0(x_0)+\beta _0(x_1)\\
\dot x_1&=&\alpha _0(x_0)-\beta _0(x_1)-\alpha _1(x_1)+\beta _1(x_2)\\
\dot x_2&=&\alpha _1(x_1)-\beta _1(x_2)-\alpha _2(x_2) -\gamma _2(x_2)+\gamma _3(x_3)\\
\dot x_3&=&\gamma _2(x_2)-\gamma _3(x_3)
\eeqn
evolving on $\R_{\geq 0}^4$, where $c\geq 0$ is a constant, and the functions 
\[
\alpha _i, \beta _i, \gamma _i : [0,\infty ) \rightarrow  [0,\infty )
\]
are all differentiable, with derivatives everywhere positive,
and so that $\alpha _i$ and $\beta _i$ are bounded, for each $i$, and $\gamma _1,\gamma _2$
are unbounded.
Furthermore, suppose that the following conditions hold:
\be{C1}
\alpha _2(\infty )+\beta _0(\infty )<\alpha _0(\infty )
\ee
\be{C2}
\alpha _0(\infty )+\beta _1(\infty )<\alpha _1(\infty )+\beta _0(\infty )
\ee
\be{C4}
\alpha _2(\infty )+\beta _1(\infty )<\alpha _1(\infty )
\ee
\be{C3}
c < \alpha _2(\infty ) \,.
\ee
Then, there is a (unique) globally asymptotically stable equilibrium for
the system.
\ets

Note that~(\ref{C1}) and~(\ref{C3}) imply also:
\be{C1p}
c+\beta _0(\infty )<\alpha _0(\infty ) \,.
\ee

\bpr
We start by noticing that solutions are defined for all $t\geq 0$.
Consider any maximal solution $x(t)=(x_0(t),x_1(t),x_2(t),x_3(t))$.
From
\be{sumP}
\frac{d}{dt} \left( x_0 + x_1 + x_2 + x_3\right) \;=\;
c - \alpha _2(x_2)
\ee
we conclude there is an estimate
$x_i(t)\leq \sum_ix_i(t)\leq \sum_ix_i(0) + tc$
and hence there are no finite escape times.
Moreover, we claim that $x(\cdot )$ is bounded.

Since the system is a strongly monotone tridiagonal system, we know that
$x_3(t)$ is {\em eventually monotone\/}.
That is, for some $T>0$, either
\be{mon1}
\dot x_3(t) \geq  0 \; \; \;\forall\,t\geq T
\ee
or
\be{mon2}
\dot x_3(t) \leq  0 \; \; \;\forall\,t\geq T \,.
\ee
Hence, $x_3(t)$ admits a limit, either finite or infinite.
Assume first that $x_3(t)\rightarrow \infty $.  Then, case~(\ref{mon2}) cannot hold,
so~(\ref{mon1}) holds.
Looking at the differential equation for $x_3$, we know that
$\gamma _2(x_2(t))-\gamma _3(x_3(t))\geq 0$ for all $t\geq T$, which means that
\[
x_2(t)\geq \gamma _2^{-1}(\gamma _3(x_3(t)))\rightarrow \infty  \,.
\]
Looking again at~(\ref{sumP}), and using that
$c-\alpha _2(\infty )<0$ (property~(\ref{C3})), we conclude that
$\frac{d}{dt} \left( x_0 + x_1 + x_2 + x_3\right)(t)<0$ for all
$t$ sufficiently large.
Thus $x_0 + x_1 + x_2 + x_3$ is bounded (and nonnegative), and this implies
that $x_2$ is bounded, a contradiction.
So $x_3$ is bounded.

Next we examine the equation for $\dot x_2$.
The two positive terms are bounded: the one involving $\alpha _1$ because
$\alpha _1$ is a bounded function, and the one involving $\gamma _3$ because $x_3$
is bounded.  Thus
\[
\dot x_2 \leq  v(t) - \alpha _2(x_2) \,,
\]
where $0\leq v(t)\leq k$ for some constant $k$.
Thus $\dot x_2(t)<0$ whenever $x_2(t)>\gamma _2^{-1}(k)$, and this proves that $x_2$
is bounded, as claimed.

Now we show that $x_0$ and $x_1$ are bounded as well.
For $x_0$, it is enough to notice that
$\dot x_0\leq c-\alpha _0(x_0)+\beta _0(\infty )$, so that
\[
x_0(t)> \alpha _0^{-1}(c+\beta _0(\infty )) \;\;\Rightarrow \;\; \dot x_0(t)<0
\]
so~(\ref{C1p}) shows that $x_0$ is bounded.
Similarly, for $x_1$ we have that
$\dot x_1\leq \alpha _0(\infty )-\beta _0(x_1)-\alpha _1(x_1)+\beta _1(\infty )$
so~(\ref{C2}) provides boundedness.

Once that boundedness has been established, if we also show that there is a
unique equilibrium then the theory of strongly monotone tridiagonal systems
(\cite{Sm, Smith}) will ensure global asymtotic stability of the equilibrium.
So we show that equilibria exist and are unique.
It is convenient to change variables are write
\[
y_0:=x_0+x_1+x_2+x_3,\;
y_1:=x_1+x_2+x_3,\;
y_2:=x_2+x_3,\;
y_3:=x_3\,.
\]
In terms of these variables, we may set $\dot y_i=0$, $i=0,2,1,3$,
so that the equilibria are precisely the solutions of:
\beqn
\alpha _2(x_2)&=&c\\
\alpha _1(x_1)&=&\alpha _2(x_2)+\beta _1(x_2)\\
\alpha _0(x_0)&=&\alpha _2(x_2)+\beta _0(x_1)\\
\gamma _3(x_3)&=&\gamma _2(x_2)\,.
\eeqn
This shows uniqueness (all the functions are strictly increasing), and
existence follows from, respectively,
(\ref{C3}), (\ref{C4}), (\ref{C1}),
and the fact that $\gamma _3$ is unbounded.
\epr

\section{Closing the Loop}

\begin{figure}[ht]
\begin{center}
\includegraphics[width=6cm,height=6cm]{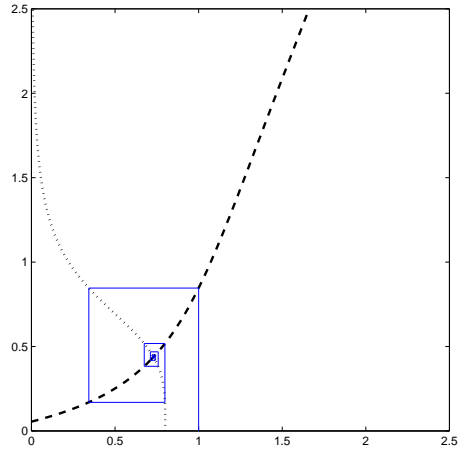}
\caption{Stability of spiderweb ($v_s=0.4$)}
\label{fig-circ4}
\end{center}
\end{figure}

Now we are ready to apply the main theorem in~\cite{monotone1}.
In order to do this, we need to plot the characteristics.
See Figure~\ref{fig-circ4} for the ``spiderweb diagram''
(the dotted and dashed curves are the characteristics)
that shows convergence of the discrete iteration described in~\cite{monotone1}
when we pick the parameter $v_s=0.4$.
The theorem implies that no oscillations can happen in that case, even under
arbitrary delays in the feedback from $P_N$ to $M$.

\begin{figure}[ht]
\begin{center}
\includegraphics[width=6cm,height=6cm]{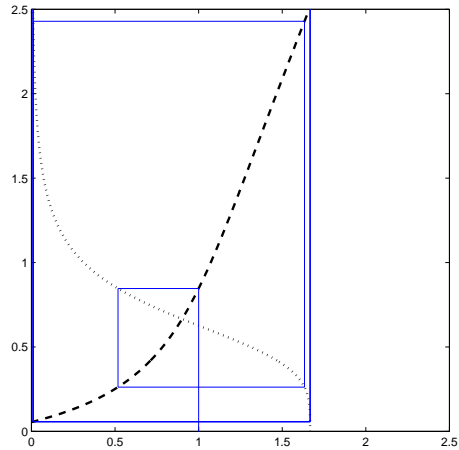}
\caption{Instability of spiderweb ($v_s=0.5$)}
\label{fig-circ5}
\end{center}
\end{figure}

\begin{figure}[ht]
\begin{center}
\includegraphics[width=6cm,height=6cm]{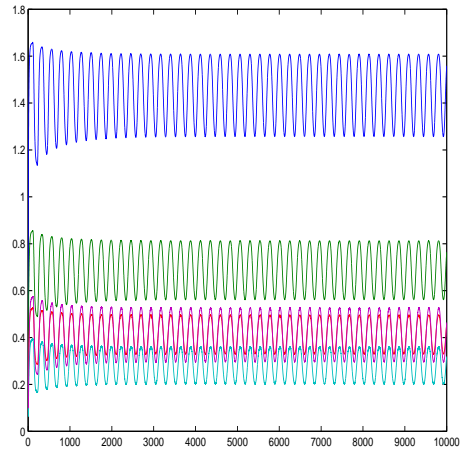}
\caption{Oscillations seen in simulations ($v_s=0.5$, delay of $100$, initial
conditions all at $0.2$), using MATLAB's dde23 package}
\label{fig-circ5-simul}
\end{center}
\end{figure}

On the other hand, for a larger value, such as $v_s=0.5$, the discrete
iteration conditions are violated;
see Figure~\ref{fig-circ5} for the ``spiderweb diagram'' that shows
divergence of the discrete iteration.
Thus, and one may
expect periodic orbits in this case.
Indeed, simulations show that, for large enough
delays, such periodic orbits arise, see Figure~\ref{fig-circ5-simul}.

\end{document}